\documentclass[pre,10pt,twocolumn,showpacs, superscriptaddress]{revtex4}
\usepackage{amsfonts}
\usepackage{amsmath}
\usepackage{amssymb}
\usepackage{graphicx}
\newtheorem{theorem}{Theorem}

\newtheorem{corollary}[theorem]{Corollary}

\newtheorem{proposition}[theorem]{Proposition}

\begin{document}
\title{On the degeneracy of $SU(3)_k$ topological phases}
\author{Stephen P. Jordan}
\affiliation{Institute for Quantum Information, California Institute of Technology}
\author{Toufik Mansour}
\affiliation{Department of Mathematics, University of Haifa}
\author{Simone Severini}
\affiliation{Department of Physics and Astronomy, University College London}

\begin{abstract}
The ground state degeneracy of an $SU(N)_k$ topological phase with $n$
quasiparticle excitations is relevant quantity for quantum
computation, condensed matter physics, and knot theory. It is an open
question to find a closed formula for this degeneracy for any $N >
2$. Here we present the problem in an explicit combinatorial way and
analyze the case $N=3$. While not finding a complete closed-form
solution, we obtain generating functions and solve some special cases.
\end{abstract}

\pacs{02.10.Ox, 31.15.xm}

\maketitle

\section{Introduction}

Topological phases have become an important object of study in both condensed
matter physics and quantum computation. Condensed matter theorists have
proposed that certain many-electron systems confined to two dimensions may
support topological phases. For example, it is widely suspected that the low
energy behavior of certain fractional quantum Hall systems is described by
the $SU(2)_k$ topological quantum field theory (TQFT) for various
$k$ \cite{Nayak}.

The world lines of $n$ particles confined to two dimensions trace out a braid
as the particles are exchanged. The induced transformation on a
$f$-fold degenerate ground space is a $f$-dimensional unitary
representation of the braid group $B_n$. If the representation is
anything other than the trivial or sign representations, the particles
are called anyons. In general, the representation could be
reducible. In this paper, we attempt to calculate the dimension of the
irreducible subspaces. These have a physical interpretation. Each
$SU(N)_k$ TQFT has a finite set of anyon types. If multiple anyons are
bound together and treated as a unit, they behave collectively as
another type from this set. If we take all $n$ anyons and fuse them
together, the type of the resulting single anyon indexes the
irreducible subspace in which the system lies. Braiding and recoupling
of anyons cannot move the system from one irreducible subspace to
another. As discussed in section \ref{setup}, for $SU(N)_k$ the irreducible
subspaces are indexed by a restricted class of Young diagrams.

In many cases, the dimensions of the irreducible subspaces grow
exponentially with $n$. To those interested in quantum computing, the
exponentially large state space, nonlocal degrees of freedom, and energy gap
all make topological phases promising candidates for fault-tolerant quantum
memories. Furthermore, by adiabatically braiding the $SU(2)_3$ anyons around
each other, one can in principle do universal quantum computation
\cite{Freedman}. This result has been extended to $SU(N)_k$ for all $N
\geq 2$ and all $k \geq 3$ other than four \cite{AharonovArad,
  WocjanYard, FLW}. The universality results all depend on
exponentially large degeneracy. The representation of the braid group
corresponding to $SU(N)_k$ also produces a topological invariant of
knots and links called the single-variable HOMFLY polynomial at the
$(N+k)^{\mathrm{th}}$ root of unity. $N=2$ yields the Jones polynomial
as a special case (see \cite{WocjanYard, AJL, ShorJordan, JordanWocjan}).

Mansour and Severini give an exact closed formula for the dimensions
of the irreducible subspaces for $SU(2)_k$ \cite{Severini}. The result
was obtained combinatorically by counting paths in the Bratteli
diagrams. The present work studies $SU(3)_k$ by similar
techniques. With these parameters the problem appears more
difficult. After introducing the problem, we write a form for the
generating function for the dimensions of the irreducible
subspaces. We then obtain explicit formulas in some special
cases. Giving a general closed formula remains an open question.

\section{Setup}
\label{setup}

\begin{figure}
[htbp]
\begin{center}
\includegraphics[width=3in]
{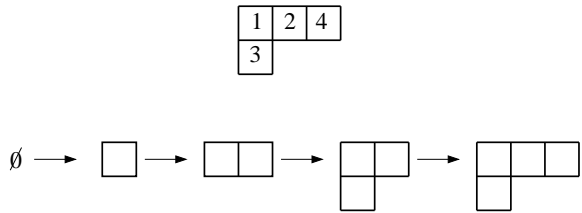}
\caption{A standard Young tableau can be interpreted as instructions for
constructing a Young diagram by adding one box at a time. After each step, the
resulting configuration is a valid Young diagram.}%
\label{fig1}
\end{center}
\end{figure}

A \emph{Young diagram} is a partition of $n$ boxes into rows, such that no row
is longer than the row above it. A \emph{standard Young tableau} is a Young
diagram in which the $n$ boxes have been numbered from $1$ to $n$. The numbers
in any column must increase downward and the numbers in any row must increase
from left to right. We can interpret these numbers as instructions for
building the Young diagram by adding one box at a time, as illustrated in Fig.
\ref{fig1}.

The condition that numbers must increase rightward and downward is equivalent
to the condition that the configuration obtained after the addition of each
box must always be a valid Young diagram. As discussed in \cite{JordanWocjan,
WocjanYard}, the ground space of a topological phase corresponding to
$SU(N)_k$ and $n$ quasiparticles separates into invariant subspaces.
Braiding and recoupling of the quasiparticles cannot move the system from one
invariant subspace to another. These subspaces correspond to the different
Young diagrams of $n$ boxes $N$ rows such that the number of boxes in the
first row minus the number of boxes in the $N^{\mathrm{th}}$ row is at most
$k$. For the $SU(N)_k$ TQFT, the dimension $f$ of the subspace with an
$n$-box Young diagram $\lambda$ is equal to the number of Young
tableaux of shape $\lambda$ such that the configuration obtained after
adding each box is not only a valid Young diagram, but also has the
property that the number of boxes in the first row minus the number of
boxes in the $N^{\mathrm{th}}$ row is at most $k$.

\begin{figure}[htbp]
\begin{center}
\includegraphics[height=3.1328in, width=3.1151in]{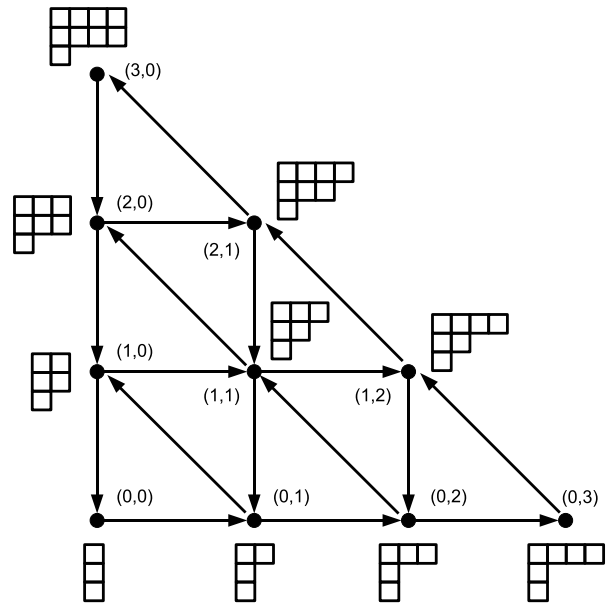}
\caption{\label{D3} The graph $D_{3}$. }
\label{fig2}
\end{center}
\end{figure}

A 3-row Young tableau can be characterized by three numbers: $n$, the
total number of boxes, $i$, the overhang of the top row over the
middle row, and $j$, the overhang of the middle row over the bottom
row. For $SU(3)_k$, $(i,j)$ is restricted to
lie within the set $V_k$ of pairs of nonnegative integers such that
$i+j \leq k$. In figure \ref{D3}, we construct a graph $D_3$ on the
vertices $V_3$. Each vertex $(i,j)$ is illustrated with an example of
a Young diagram with the corresponding set of overhangs. A directed
edge from $(i,j)$ to $(i',j')$ is included if one can go from a Young
diagram with overhangs $(i,j)$ to a Young diagram with overhangs
$(i',j')$ by adding one box. The $n$-box Young tableaux of shape
$\lambda$ allowed for $SU(3)_k$ correspond bijectively to the paths on
the graph $D_k$ starting from $(0,0)$ and ending after $n$ steps on
the vertex corresponding to $\lambda$. More precisely, a \emph{path}
of length $n$ in $D_{k}$ is a sequence of vertices
$v_{0},v_{1},\ldots,v_{n}\in D_{k}$ and edges $\left(
v_{0},v_{1}\right) ,\left(  v_{1},v_{2}\right)  ,\ldots,\left(
v_{n-2},v_{n-1}\right) ,\left(  v_{n-1},v_{n}\right)  \in
D_{k}$. A path can contain more than a single occurrence of the same vertex.

Let $f_{i,j}(n,k)$ be the number of paths on $D_k$ starting from
$(0,0)$, and ending on $(i,j)$ after $n$ steps. $f_{i,j}(n,k)$ is
equal to the dimension of the invariant subspace of $n$ $SU(3)_k$
anyons whose Young diagram has overhangs $(i,j)$. The remainder of
this paper is devoted to analyzing $f_{i,j}(n,k)$.

The set of directed edges $A_k$ in the graph $D_k$ can be described formally
by the following constraints. $\left(  \left(  a,b\right)  ,\left(
c,d\right)  \right)  \in A_{k}$ only in the following cases:

\begin{itemize}
\item if $a=c$ and $b\neq d$ then $d=b+1$;
\item if $a\neq c$ and $b=d$ then $c=a-1$;
\item if $a\neq c$ and $b\neq d$ then $c=a+1$ and $d=b-1$.
\end{itemize}

The adjacency matrix of a graph is a
matrix in which the $ij^{\mathrm{th}}$ entry is $1$ if there is an edge $\left(
i,j\right)  $, otherwise it is $0$. The number of paths of length $n$ from
vertex $i$ to vertex $j$ equals the $ij^{\mathrm{th}}$ entry in the
$n^{\mathrm{th}}$ power of the adjacency matrix.

\section{Results}

Table~\ref{tff} gives $f_{0,0}(n,k)$ for all
$3|n\leq27$. $f_{0,0}(n,k) = 0$ for all $l$ not divisible by
three. Notice that the diagonal entries of the table are the
$3$-dimensional Catalan numbers $\frac{2 n!}{(n/3)!(n/3+1)!
  (n/3+2)!}$. The corresponding generating function is denoted by
\begin{equation}
\label{generate}
F_{i,j}(t;k)=\sum_{n\geq0}f_{i,j}(n,k)t^{n}.
\end{equation}

\begin{table}[htbp]
\caption{Values of $f_{0,0}(n,k)$ for $k \leq 8$ and $n \leq 9$}
\begin{tabular}
[c]{c|cccccccccc}
$k\backslash n$ & $0$ & $3$ & $6$ & $9$ & $12$ & $15$ & $18$ & $21$ & $24$ &
$27$\\\hline\hline
$1$ & $1$ & $1$ & $1$ & $1$ & $1$ & $1$ & $1$ & $1$ & $1$ & $1$\\
$2$ & $1$ & $1$ & $5$ & $21$ & $89$ & $377$ & $1597$ & $6765$ & $28657$ &
$121393$\\
$3$ & $1$ & $1$ & $5$ & $42$ & $341$ & $2731$ & $21846$ & $174763$ & $1398101$
& $11184810$\\
$4$ & $1$ & $1$ & $5$ & $42$ & $462$ & $5278$ & $60181$ & $683962$ & $7763097$
& $88079511$\\
$5$ & $1$ & $1$ & $5$ & $42$ & $462$ & $6006$ & $83028$ & $1166677$ &
$16440171$ & $231612211$\\
$6$ & $1$ & $1$ & $5$ & $42$ & $462$ & $6006$ & $87516$ & $1357569$ &
$21669957$ & $349920000$\\
$7$ & $1$ & $1$ & $5$ & $42$ & $462$ & $6006$ & $87516$ & $1385670$ &
$23193775$ & $401389561$\\
$8$ & $1$ & $1$ & $5$ & $42$ & $462$ & $6006$ & $87516$ & $1385670$ &
$23371634$ & $413180625$\\
\end{tabular}
\label{tff}
\end{table}

From the definitions, we can state that
\begin{align}
&  F_{i,j}(t;k)-\delta_{i=j=0}\nonumber\label{eqrec}\\
&  \quad=t(F_{i+1,j}(t;k)+F_{i-1,j+1}(t;k)+F_{i,j-1}(t;k)),
\end{align}
for all $0\leq i+j\leq k$, with the initial condition $F_{i,j}(t;k)=0$, for
all $i+j>k$, $i<0$ or $j<0$. For example, if $k=1$ then the above recurrence
relation yields $F_{0,0}=1+tF_{1,0}$, $F_{1,0}=tF_{0,1}$ and
$F_{0,1}=tF_{0,0}$. These imply $F_{0,0}(t;1) = 1/\left(
1-t^{3}\right)  $. Taylor expanding $F_{0,0}(t;1)$ and using
Eq. \ref{generate} reproduces the first row of table \ref{tff}. In order
to write a system of equations on the variables $F_{i,j}(t;k)$, let
$J_{p,q;s}$ be the $p \times q$ matrix $J_{p,q;s}(i,j)$, where
\[
J_{p,q;s}(i,j)=\left\{
\begin{array}
[c]{ll}%
1, & j-i=s;\\
0, & \mbox{otherwise}.
\end{array}
\right.
\]
Moreover, let us define the following matrices.
\begin{eqnarray*}
A_k & = & J_{k,k;0}-tJ_{k,k,-1} \\
E_p & = & -tJ_{p,p-1;0} \\
E^{\prime}_p & = & -tJ_{p,p+1;1}
\end{eqnarray*}
Let $x$ be a vector of $\binom{k+2}{2}$ coordinates defined by
\begin{align*}
&  (F_{0,0}(t;k),F_{0,1}(t;k),\ldots,F_{0,k}(t;k),\\
&  \qquad\qquad\,\,\,\,F_{1,0}(t;k),\ldots,F_{1,k-1}(t;k),\\
&  \qquad\qquad\,\,\,\,\qquad\qquad\,\,\,\,\ddots\\
&  \qquad\qquad\qquad\qquad\qquad\,\,\,\,\,\,F_{k,0}(t;k)),
\end{align*}
that is,
\[
x_{i(2k-i+3)/2+j+1}=F_{i,j}(t;k),
\]
for all $0\leq i+j\leq k$. In addition, let%
\[
F_{k}=\left(
\begin{array}
[c]{llllll}%
A_{k+1} & E_{k+1} &  &  &  & \\
E^{\prime}_{k} & A_{k} & E_{k} &  &  & \\
& E^{\prime}_{k-1} & A_{k-1} & E_{k-1} &  & \\
&  & \ddots & \ddots & \ddots & \\
&  &  & E^{\prime}_{2} & A_{2} & E_{2}\\
&  &  &  & E^{\prime}_{1} & A_{1}%
\end{array}
\right)  .
\]
Hence, rewriting Eq. \eqref{eqrec} in matrix form, we obtain the
following result.

\begin{proposition}
\label{th1}The generating functions
\[
x_{i(2k-i+3)/2+j+1}=F_{i,j}(t;k)
\]
with $0\leq i+j\leq k$ satisfy
\[
F_{k}\cdot(x_{1},\ldots,x_{\binom{k+2}{2}})^{T}=(1,0,\ldots,0)^{T}.
\]
\end{proposition}

Proposition~\ref{th1} for $k=1$ gives
\[
F_{1}^{-1}=\frac{1}{1-t^{3}}\left(
\begin{array}
[c]{lll}%
1 & t^{2} & t\\
t & 1 & t^{2}\\
t^{2} & t & 1
\end{array} \right) ,
\]
which implies:
\begin{eqnarray*}
F_{0,0}(t;1) & = & \frac{1}{1-t^{3}}, \\
F_{0,1}(t;1) & = & \frac{t}{1-t^{3}},\\
F_{1,0}(t;1) & = & \frac{t^{2}}{1-t^{3}}.
\end{eqnarray*}
Proposition~\ref{th1} for $k=2$ gives
\[
F_{2}^{-1}=\frac{1}{1-4t^{3}-t^{6}}G_{2},
\]
where
\[
G_{2}=\left(
\begin{array}
[c]{cccccc}%
1-3t^{3} & tz & 2t^{4} & ty & 2t^{3} & t^{2}y\\
ty & y & tz & 2t^{2} & z & 2t^{3}\\
t^{2}y & ty & 1-3t^{3} & 2t^{3} & tz & 2t^{4}\\
tz & z & 2t^{3} & y & 2t^{2} & ty\\
2t^{3} & 2t^{2} & ty & z & y & tz\\
2t^{4} & 2t^{3} & t^{2}y & tz & ty & 1-3t^{3}%
\end{array}
\right)  ,
\]
with $y=1-t^{3}$ and $z=t(1+t^{3})$. Therefore
\[%
\begin{array}
[c]{ll}%
F_{0,0}(t;2)=\dfrac{1-3t^{3}}{1-4t^{3}-t^{6}}, & F_{0,1}(t;2)=\dfrac
{t(1-t^{3})}{1-4t^{3}-t^{6}}\\
F_{0,2}(t;2)=\dfrac{t^{2}(1-t^{3})}{1-4t^{3}-t^{6}}, & F_{1,0}(t;2)=\dfrac
{t(1+t^{3})}{1-4t^{3}-t^{6}}\\
F_{1,1}(t;2)=\dfrac{2t^{3}}{1-4t^{3}-t^{6}}, & F_{2,0}(t;2)=\dfrac{2t^{4}%
}{1-4t^{3}-t^{6}}.
\end{array}
\]
This implies
\[
f_{0,0}(n;2) = \left\{ \begin{array}{cl}
\mathrm{Fib}_{n-1} & \textrm{if $3|n$} \\
0                  & \textrm{otherwise}
\end{array} \right.
\]
where $\mathrm{Fib}_n$ is the $n^{\mathrm{th}}$ Fibonacci number.
Applying proposition~\ref{th1} for $k=1,2,3,4,5$, we have the next corollary.

\begin{corollary}
The generating function for the number of paths of length $n$ from $(0,0)$
to $(0,0)$ in $D_{k}$, $k=1,2,3,4$, is given by

\begin{itemize}
\item $F_{0,0}(t;1)=\frac{1}{1-t^{3}}$,

\item $F_{0,0}(t;2)=\frac{1-3t^{3}}{1-4t^{3}-t^{6}}$,

\item $F_{0,0}(t;3)=\frac{1-8t^{3}+5t^{6}-2t^{9}}{1-9t^{3}+9t^{6}-8t^{9}}$,

\item $F_{0,0}(t;4)=\frac{1-15t^{3}+48t^{6}-46t^{9}-19t^{12}}{1-16t^{3}%
+59t^{6}-67t^{9}-37t^{12}+8t^{15}}$.

\end{itemize}
\end{corollary}

From the definition of determinant, we see that $\det(F_k)$ has the
following properties.
\begin{itemize}
\item $\det(F_{k})$ is a polynomial of degree $d_{k}$, where $d_{3k-1}%
=\frac{3k(3k+1)}{2}$, $d_{3k}=\frac{9k(k+1)}{2}$ and $d_{3k+1}=\frac
{3(k+1)(3k+2)}{2}$.

\item $\det(F_{k})$ is given by $1-k^{2}t^{3}+t^{6}p_{k}(x)$, where $p_{k}(x)$
is a polynomial.
\end{itemize}
This implies that $\rho_k$, the smallest positive root of the polynomial
$\det(F_{k})$, is approximated by $k^{-2/3}$. If $N_{k,n}(i,j)$ is the
number of paths of length $n$ from $(0,0)$ to $(i,j)$ in $D_{k}$ then
\[
\lim_{n\rightarrow\infty}(N_{k,n}(i,j))^{1/n}=\frac{1}{\rho_{k}}.
\]\vspace{-10pt}

\begin{table}[!htbp]
\caption{Values of $\det(F_{k})$ for $k=1,2,\ldots,8$}
\begin{tabular}
[c]{c|l}
$k$ & $\det(F_{k})$\\\hline\hline
$1$ & $1-t^{3}$\\\hline
$2$ & $1-4t^{3}-t^{6}$\\\hline
$3$ & $1-9t^{3}+9t^{6}-8t^{9}$\\\hline
$4$ & $1-16t^{3}+59t^{6}-67t^{9}-37t^{12}+8t^{15}$\\\hline
$5$ & $1-25t^{3}+191t^{6}-559t^{9}+531t^{12}-507t^{15}+341t^{18}$\\
& $+27t^{21}$\\\hline
$6$ & $1-36t^{3}+459t^{6}-2655t^{9}+7290t^{12}-9801t^{15}$\\
& $+3429t^{18}+6075t^{21}-1458t^{24}+729t^{27}$\\\hline
$7$ & $1-49t^{3}+929t^{6}-8865t^{9}+46315t^{12}-136058t^{15}$\\
& $+219202t^{18}-198802t^{21}+189535t^{24}-152085t^{27}$\\
& $+62341t^{30}+20851t^{33}-1331t^{36}$\\\hline
$8$ & $1-64t^{3}+1679t^{6}-23699t^{9}+198636t^{12}-1031272t^{15}$\\
& $+3360456t^{18}-6855112t^{21}+8542281t^{24}-5062167t^{27}$\\
& $-1959023t^{30}+4912958t^{33}-1335971t^{36}+1092507t^{39}$\\
& $-375746t^{42}-12167t^{45}$\\\hline
\end{tabular}
\end{table}


For large $n$ the number of paths of length $n$ in $D_k$ from $(0,0)$
to any other vertex scales as $\lambda_k^n$ where $\lambda_k = 1/\rho_k$
is the largest eigenvalue of the adjacency matrix of $D_k$. This
quantity coincides with the ``total quantum dimension'' of the TQFT
$SU(N)_k$, which is given\cite{WocjanYard} by
\begin{equation}
\lambda_k = \frac{\sin(\pi N/(N+k))}{\sin(\pi/(N+k))}.
\end{equation}

In the case that $k \geq n$ the restriction that total overhang $i+j$
is at most $k$ becomes irrelevant. In this case our problem reduces
to counting ordinary Young tableaux without any special
restrictions. The solution to this problem can be derived from the
hook length formula as described in proposition \ref{hooks}.

\begin{proposition}
\label{hooks}
The number of paths of length $n$ from vertex $(0,0)$ to vertex $(i,j)$ in
$D_{k \geq n}$ is
\begin{equation}
\frac{(i+1)(j+2)(j-i+1)n!}{(\frac{n-i+2j+6}{3})!(\frac{n+2i-j+3}%
{3})!(\frac{n-i-j}{3})!}, \label{numb}%
\end{equation}
where $(n-i+2j)/3$ is a positive integer. Otherwise, there is no such path.
\end{proposition}

Note that proposition \ref{th1} implies that for any $k$, the generating function
$F_{i,j}(t;k)$ is a rational function on $t$. Finding general explicit formulas for
the generating function $F_{i,j}(t;k)$ and the dimension
$f_{i,j}(n,k)$ remain open problems. Additional physically motivated
open problems include computing the total dimension $\sum_{i,j}
f_{i,j}(n;k)$ and investigating the $N \to \infty$ limit, which one
might expect to exhibit semiclassical behavior.

\emph{Acknowledgements: } S.J. acknowledges support from the Sherman
Fairchild Foundation and NSF grant PHY-0803371.


\begin{thebibliography}{9}                                                                                                %
\bibitem {Nayak}C. Nayak, S. H. Simon, A. Stern, M. Freedman, and S. Das
Sarma, Non-Abelian anyons and topological quantum computation, \emph{Rev. Mod.
Phys.} \textbf{80}, 1083 (2008). arXiv:0707.1889v2 [cond-mat.str-el]

\bibitem {Freedman}M. Freedman, M. Larsen, and Z. Wang, A modular functor
which is universal for quantum computation, \emph{Comm. Math. Phys.}
\textbf{227}, no. 3, 605-622 (2002). arXiv:quant-ph/0001108v2

\bibitem {FLW}M. Freedman, M. Larsen, and Z. Wang, The two-eigenvalue
  problem and density of Jones representation of braid
  groups. \emph{Comm. Math. Phys.} \textbf{228}, no. 1,
  177-199. (2002) arXiv:math.GT/0103200.

\bibitem {AharonovArad}D. Aharonov and I. Arad, The BQP-hardness of
approximating the Jones Polynomial, 2006. arXiv:quant-ph/0605181v2

\bibitem {WocjanYard}P. Wocjan and J. Yard, The Jones polynomial: quantum
algorithms and applications in quantum complexity theory, \emph{Quantum
Information and Computation}, \textbf{8}:1-2 (2008), 147-180. arXiv:quant-ph/0603069

\bibitem {AJL}D. Aharonov, V. Jones and Z. Landau, A Polynomial Quantum
Algorithm for Approximating the Jones Polynomial, \emph{Proceedings of the
38th ACM Symposium on Theory of Computing}, 2006. arXiv:quant-ph/0511096

\bibitem {JordanWocjan}S. P. Jordan and P. Wocjan, Estimating Jones and HOMFLY
polynomials with One Clean Qubit, \emph{Quantum Information and Computation},
\textbf{9}:3-4 (2009), 264-289. arXiv:0807.4688

\bibitem {Severini}T. Mansour and S. Severini, Counting paths in Bratelli
diagrams for $SU\left(  2\right)  _{k}$, \emph{Europhysics Letters} \textbf{86}(2009), 33001. arXiv:0806.4809

\bibitem {ShorJordan}P. Shor and S. Jordan, Estimating Jones polynomials is a
complete problem for one clean qubit, \emph{Quantum Information and
Computation}, \textbf{8}:8-9 (2008), 181-214. arXiv:0707.2831
\end{thebibliography}
\end{document}